\documentclass[12pt]{amsart} \usepackage{amscd}

\newtheorem{theorem}{Theorem}
\newtheorem{lemma}[theorem]{Lemma}
\newtheorem{corollary}{Corollary}

\theoremstyle{definition}

\theoremstyle{remark} \newtheorem*{remark}{Remark}

\usepackage{amsfonts}
\newcommand{\field}[1]{\mathbb{#1}}          \newcommand{\Q}{\field{Q}}
\newcommand{\R}{\field{R}}                   \newcommand{\Z}{\field{Z}}
\newcommand{\C}{\field{C}}                   \newcommand{\A}{\field{A}}

\newcommand{\fg}{\mathfrak     g}     \newcommand{\fp}{\mathfrak    p}
\newcommand{\fk}{\mathfrak     k}     \newcommand{\fh}{\mathfrak    h}
\newcommand{\fm}{\mathfrak m}

\newcommand{\bs}{\backslash} \newcommand{\ra}{\rightarrow}

\begin{document}

\title[Construction of  Some Generalised Modular Symbols]{Construction
of Some Generalised Modular Symbols}

\author{B.Speh and T. N. Venkataramana}

\address{Department   of  Mathematics,   310   Malott  Hall,   Cornell
University, Ithaca N.Y. 14853-4201, U.S.A}

\address{School   of  Mathematics,   Tata  Institute   of  Fundamental
Research, Homi Bhabha Road, Bombay - 400 005, INDIA.}

\email{speh@math.cornell.edu, venky@math.tifr.res.in}

\subjclass{Primary 11F75;  Secondary 11F80\\ {\bf  B.Speh}, Department
of   Mathematics,  310   Malott  Hall,   Cornell   University,  Ithaca
N.Y.14853-4201,   U.S.A\\    speh@math.cornell.edu\\   {\bf   T.    N.
Venkataramana}, School  of Mathematics, Tata  Institute of Fundamental
Research,   Homi   Bhabha   Road,   Bombay   -   400   005,   INDIA.\\
venky@math.tifr.res.in}

\date{}

\begin{abstract}
We give a  criterion for the non-vanishing of  certain modular
symbols on  a locally symmetric  manifold. The  criterion is  in
terms  of the non-vanishing of  some cohomology classes  on the
compact dual  of the locally  symmetric  manifold.   Using   this,
we  construct nonzero  modular symbols  for $SL_{2n}$ over an
imaginary quadratic extension of $\Q$, which represent ghost
classes. We  also construct nonzero modular symbols in certain
non-compact Shimura varieties and give an example of a modular
symbol that generates an infinite dimensional module under the
action of the Hecke algebra.
\end{abstract}

\maketitle

\section{Introduction}

Let  $G$  and  $H$  be  semi-simple algebraic  groups  over  $\Q$
and $f:~H\rightarrow G$  a morphism  with finite kernel  of
$\Q$-algebraic groups. Let $K_H$  be a maximal compact subgroup of
the group $H(\R)$ of real points of $H$ and $K\supset f(K_H)$ a
maximal compact subgroup of  $G(\R)$. We  then get  a map
$Y\rightarrow X$  of  the associated symmetric spaces
$Y=H(\R)/K_H$ and $X=G(\R)/K$. Suppose  that $X$ and $Y$ are of
dimension $D$ and $d$ respectively.

\medskip
Assume  that  $\Gamma  \subset  G(\Q)$ is  a  torsion-free
congruence subgroup. Then  $S(\Gamma )=\Gamma  \bs X$ is  a
manifold  with finite volume under the  $G(\R)$ invariant metric
(and volume  form ) on $X$. Write $\Gamma  \cap H$ for
$f^{-1}(\Gamma )$  and $S_H(\Gamma )=\Gamma \cap H\bs  Y$. We have
an immersion   $S_H(\Gamma )\rightarrow  S(\Gamma )$ induced by
the map $f$. This identifies $S_H(\Gamma )$ with a closed subspace
of $S(\Gamma )$. This follows as in \cite{A} 2.7. For ease of
notation, we will assume that the groups $G(\R)$ and $H(\R)$ are
connected. Then, both the symmetric spaces $X$ and $Y$ are
orientable with an orientation fixed by $G(\R)$ and $H(\R)$
respectively.   This   in turn makes $S(\Gamma)$ and $S_H(\Gamma
)$ orientable.

\medskip

Now,  compactly supported  cohomology  classes in $H^d_c(S(\Gamma
))$  may be pulled back to $S_H(\Gamma )$ and  integrated  on
$S_H(\Gamma )$.  Since $H^d_c(S_H(\Gamma ))=\C$ integration on
$S_H(\Gamma )$ defines a linear  form on $H^d_c(S(\Gamma  ))$
which  may  be thought  of  -by Poincar\'{e} duality-  as an
element $[S_H(\Gamma )]$ of $H^{D-d}(S(\Gamma ))$. It is called
the {\bf generalized modular symbol} corresponding to $H$.

\medskip

In the arithmetic of modular curves and classical automorphic
forms modular symbols have served as an indispensable tool linking
geometry and arithmetic. "Period integrals" of Eisenstein classes
or cuspidal cohomology classes over compact modular symbols have
been used by G.Harder to obtain information about special values
of L-functions \cite{H1}, \cite{H2}. In 1990 A. Ash and A.Borel
showed that the Levi factors of parabolic subgroups define nonzero
modular symbols \cite{A-B}, \cite{R-S}. Later Ash, Ginzburg and
Rallis give 6 families of pairs (G,H) where they can show that
that  any cuspidal cohomology class for $\Gamma $ over a
generalized modular symbol corresponding to H. \cite{A-G-R}. One
such pair is $G=Sp_n$ and $H=Sp_m \times Sp_k$ which we also
consider in theorem \ref{symplectic}.

\medskip

 In this paper  we give a criterion for the  non-vanishing of a
modular symbol $[S_H(\Gamma  )]$.  The  criterion is in terms of
the compact duals  spaces $\hat{X}$ and $\hat {Y}$ of the
symmetric spaces $X$ and $Y$. We  recall the construction of
$\hat{X}$ and $\hat {Y}$; there is a Cartan decomposition of the
Lie algebra $\fg _0$ of $G(\R)$ with respect to the maximal
compact $K$, which leaves the Lie algebra $\fh _0$ of $H(\R)$
stable. Write $\fg _0=\fk _0\oplus \fp _0$ and $\fh _0=\fk _0\cap
\fh _0\oplus \fp _0\cap \fh _0$ as the Cartan decompositions. Let
$G_u$ be the compact  subgroup of $G(\C)$ with Lie  algebra $\fk
_0\oplus  i  \fp  _0$. Then $\widehat {X}=G_u/K$.   The compact
dual $\widehat {Y}$ is defined analogously. One has an embedding
of $\widehat {Y}$ in $\widehat {X}$, and  the fundamental class
$[\widehat {Y}]$ of $\widehat {Y}$  in $\widehat {X}$  is   a
cohomology   class in $H^{D-d}(\widehat {X})$.

\medskip
The {\bf Borel  map} $j$ from the cohomology  $H^*(\widehat {X})$
into the cohomology  $H^*(S(\Gamma ))$ is  defined as follows:  we
identify $H^*(\widehat  {X})$ with $H^*(\fg,K,\C)$  and
$H^*(S(\Gamma  ))$ with $H^*(\fg, K,  C^{\infty}(\Gamma \bs
G(\R)))$  where $\fg=\fg _0\otimes \C$; then  the map  $j$ is
induced by the  inclusion of  the constant functions $\C$  in the
space $C^{\infty}(\Gamma \bs  G(\R))$ of smooth functions on
$\Gamma \bs G(\R)$.

\medskip
With notation as above, we prove the following theorem in section
2.
\begin{theorem}\label{nonvan}
Suppose  that   $G$  is  a   simply  connected  group  which   has  no
$\R$-anisotropic   factors  defined  over   $\Q$.   Then,   the  class
$j([\widehat {Y}])$ is a linear combination of Hecke translates of the
generalised modular  symbol $[S_H(\Gamma ')]\in  H^{D-d}(S(\Gamma '))$
for some congruence  subgroup $\Gamma '$ of $\Gamma  $. In particular,
if  $j([\widehat {Y}])\neq  0$, then  the modular  symbol $[S_H(\Gamma
')]$ does not vanish.
\end{theorem}

The proof of this result relies on the work of J. Franke \cite{F}.

\medskip In section 3  we deduce, from Theorem \ref{nonvan} and
some computations in the cohomology  of classical compact
symmetric spaces, the following theorem.

\begin{theorem}\label{examples}
If $E$ is a  totally imaginary number field, $G=R_{E/\Q}(SL_{2n})$ and
$H=R_{E/\Q}(Sp_{2n})$, then  the modular symbol  $[S_H(\Gamma )]$ does
not vanish for some congruence subgroup $\Gamma $.

If $E$  is a totally imaginary  quadratic extension of  a totally real
number  field  $F$,  $G=R_{E/\Q}SL_{2n+1}$ and  $H=R_{F/\Q}SL_{2n+1}$,
then  the modular  symbol $[S_H(\Gamma  )]$ does  not vanish  for some
congruence subgroup $\Gamma $.
\end{theorem}
In the above, $R$ denotes the (Weil) restriction of scalars.

\medskip A cohomology class was called a {\it ghost class } by A.Borel
if  it restricts  trivially to  each boundary  component of  the Borel
Serre compactification,but its restriction to the full boundary is not
zero. The first example of  ghost class was constructed by G.Harder in
the  cohomology   of  $GL_3$  over  totally   imaginary  fields  using
Eisenstein  classes.   Later  a  whole  family of  ghost  classes  was
constructed  by  J. Franke  using  cohomology  classes represented  by
invariant forms. We show

\begin{corollary}\label{ghost}
The $G(\A_f)$-span of the  modular symbols $[S_H(\Gamma )]$ in theorem
\ref{examples} contains ghost classes.
\end{corollary}

Theorem  \ref{nonvan}  is  especially  useful  in the  case  when  the
symmetric spaces $X$  and $Y$ are of Hermitian  type and the embedding
$Y\ra X$ is holomorphic.  In particular, we prove
\begin{theorem}\label{symplectic}
Let $G=Sp_{2g}$ be the split  symplectic group over $\Q$ and let $H=\prod
Sp_{2g_i}\subset  Sp_{2g}$ with  $\sum g_i=g$,  be the  natural inclusion.
Then, the modular symbol $[S_H(\Gamma  )]$ is non-zero, for a suitable
congruence subgroup $\Gamma $.
\end{theorem}

Analogously, for the unitary group, we prove
\begin{theorem}\label{unitary}
Suppose  $q\geq p\geq  1$  are  integers, and  let  $G=U(p,q)$ be  the
unitary group in $p+q$ variables.

(1) If $p_i$ and  $q_i$ are integers such that  $\sum p_i=p$ and $\sum
    q_i\leq  q$,  and $H=\prod  U(p_i,q_i)$  then  the modular  symbol
    $[S_H(\Gamma )]$ is non-zero, for some congruence subgroup $\Gamma
    $ of $G$.

(2) Under the natural embedding of $H=Sp_{2g}$ in $G=U(g,g)$, the modular
symbol  $[S_H(\Gamma  )]$ is  non-zero  for  some congruence  subgroup
$\Gamma $.
\end{theorem}

In some instances, one can even  determine if the space of the span of
$G(\A_f)$-translates  of  the   modular  symbol  $[S_H(\Gamma  )]$  is
infinite dimensional.
\begin{theorem}\label{infinite}
Suppose that $G=U(1,q)$, and  $H=U(1,r)$ with $r=q-2$ or $r=q-1$. Then
there  exists   a  congruence  subgroup   $\Gamma  $  such   that  the
$G(\A_f)$-translates  of  the   modular  symbol  $[S_H(\Gamma  )]$  is
infinite dimensional.
\end{theorem}

\section{some results of Franke and proof of Theorem \ref{nonvan}}

\subsection{Notation} Let $\Gamma $ be as in Theorem \ref{nonvan}. Let
$K_f$ be the  closure of $\Gamma $ in  $G(\A_f)$, where $\A_f$ denotes
the ring  of finite adeles over  $\Q$.  Put $H^*(S_G)$  for the direct
limit  $lim  H^*(S(\Gamma ))$  as  $\Gamma  $  varies over  congruence
subgroups of  $G(\Q)$ (if $\Gamma '\subset  \Gamma $, then  there is a
natural inclusion $H^*(S(\Gamma ))$  in $H^*(S(\Gamma '))$ which gives
us a  direct system of  finite dimensional complex vector  spaces, and
the direct limit is with  respect to these inclusions). On this direct
limit,  the group $G(\A_f)$  operates, and  if $K$  is the  closure of
$\Gamma  $  in $G(\A_f)$  then,  the  space  of $K_f$  -invariants  in
$H^*(S_G)$  is   exactly  $H^*(S(\Gamma  ))$  (we   are  using  strong
approximation  here  which  is  guaranteed under  the  assumptions  of
Theorem \ref{nonvan}).

\medskip  Let  $K_0$ be  a  {\bf  good  maximal compact}  subgroup  of
$G(\A_f)$. Put  $\Gamma _0=G(\Q)\cap K_0$.  We fix  a subgroup $\Gamma
'$ of finite  index in $\Gamma \cap \Gamma _0$.  It  is clear from the
definition of the modular  symbol $\xi _{\Gamma }=[S_H(\Gamma )]$ that
the  sum of  translates  of $\xi  _{\Gamma  '}$ over  a  set of  coset
representatives of $\Gamma /\Gamma  '$ (resp. $\Gamma _0/\Gamma '$) is
a non-zero multiple  of $\xi _{\Gamma }$ (resp.   $\xi _{\Gamma _0}$).
Therefore, if we show that $j([\widehat {Y}])$ is a linear combination
of $G(\A_f)$-translates of $S_H(\Gamma _0)$ then we have proved Therem
\ref{nonvan}.

Set  ${\mathcal H}_0$  to be  the  space of  complex valued  compactly
supported  $K_0$-bi-invariant  functions on  the  finite adelic  group
$G(\A_f)$. This  is an algebra (Hecke algebra  corresponding to $K_0$)
under  convolutions and  acts  on the  cohomology group  $H^*(S(\Gamma
_0))=H^*(S_G)^{K_0}$.   Let $\C$  denote the  trivial  one dimensional
$G(\A_f)$-module. On  this module,  ${\mathcal H}_0$ operates,  and we
get a  homomorphism $\chi :{\mathcal  H}_0\rightarrow \C$. Let  $\fm =
\fm _{\chi}$  denote the kernel of  the map $\chi$. This  is a maximal
ideal in ${\mathfrak H}_0$. Denote by $H^*(S_G)^{K_0}_{\fm}$ the space
of  vectors in  $H^*(S_G)^{K_0}$ which  are annihilated  by  some {\bf
power} of the ideal $\fm$. Clearly,  this space is a direct summand as
a Hecke module.

However, the main theorem (Theorem  1 of \cite{F}) of \cite{F} asserts
that for the good maximal compact subgroup $K_0$, this space coincides
with the space of vectors annihilated  by the {\bf first} power of the
maximal ideal  $\fm$ and that  it is also  isomorphic to the  space of
co-invariants  for $G(\A_f)$  of  the module  $H^*(S_G)$ .   Moreover,
according  to   \cite{F},  the  latter  is   naturally  isomorphic  to
$H^*({\mathcal  U})^{K\cap  P(\R)}$  where  $P$  is  a  fixed  minimal
parabolic  $\Q$-subgroup of $G$  and ${\mathcal  U}$ is  the following
open set  in the compact dual $\widehat  {X}$: \[{\mathcal U}=\widehat
{X}-\cup \widehat {X_L}\]  where the union is over  Levi subgroups $L$
of  parabolic $\Q$-subgroups of  $G$ containing  the minimal  one $P$;
$\widehat  {X_L}$ denotes  the  compact dual  of  the symmetric  space
associated  to  the  group  $L^0\subset  L(\R)$. Here,  $L^0$  is  the
subgroup of $L$, which is  the intersection of the kernels of rational
characters on  $L$.  The above isomorphism takes  the restriction $res
([x])$ to  $H^*({\mathcal U})$ of a class  $[x]\in H^*(\widehat {X})$,
into the vector $j([x])$. \\

Now  the  space $V=H^*(S_G)^{K_0}$  is  a  finite dimensional  complex
vector space which is a module over the ring ${\mathcal H}_0$. Let $J$
be  the  kernel  of  the   map  ${\mathcal  H}_0\ra  End(V)$  and  let
$R={\mathcal  H}_0/J$  be  the   quotient  ring  (which  is  a  finite
dimensional algebra over $\C$). Then,  $\fm$ is the inverse image of a
maximal ideal $\fm  _R$ of $R$ under the quotient  map. Let $I$ denote
the annihilator  in $R$ of  the modular symbol $[S_H(\Gamma  _0)]$. We
may write  $I=\fm _R^kA$ where  $A$ is an  ideal coprime to  $\fm _R$.
The span $W_{[S_H(\Gamma _0)]}$  of ${\mathcal H}_0$-translates of the
class $[S_H(\Gamma _0)]$  is isomorphic as an $\mathcal  H_ 0$ module,
to $R/I=R/\fm _R^k\oplus R/A$. By  Franke's Theorem, if $k\geq 1$ then
$k=1$. Consequently,  if the  projection $I([S_H(\Gamma _0)])$  of the
modular symbol  $[S_H(\Gamma _0)]$ to  the space of  ``invariants'' is
non-zero,   then  the   projection  $I([S_H(\Gamma   _0)])$   lies  in
$W_{[S_H(\Gamma _0)]}$ and  hence it is a linear  combination of Hecke
translates of $[S_H(\Gamma _0)]$.

We will  complete the  proof of Theorem  \ref{nonvan} by  showing that
this  projection  $I([S_H(\Gamma  _0)])$  is  a  nonzero  multiple  of
$j([\widehat {Y}])$.\\

The space of  $\fm $ invariants is a direct  summand of the ${\mathcal
H}_0$-module  $H^*(S(\Gamma _0))$. We  identify it  as a  vector space
using the map j with $H^*({\mathcal U})$. Hence under Poincare duality
we  can  identify $H^*_c({\mathcal  U})$  with  the  space of  $\fm  $
invariants  in  $H^*_c(S(\Gamma  _0))$.   Note  that  since  the  ring
${\mathcal  H}_0$ is  equipped  with a  nice  involution induced  from
$g\mapsto g^{-1}$ on  $G(\A_f)$, we can consider these  spaces as dual
modules of ${\mathcal H}_0$-modules.

For a  given class $[w]$  in $j(H^d_c({\mathcal U}))$ we  consider the
map \[W_{[S_H(\Gamma _0)]} \rightarrow H^D_c(S(\Gamma _0))\] defined $
[v] \rightarrow [v]  \wedge [w] $ as a  linear form on $W_{[S_H(\Gamma
_0)]}$.   Since [w] is  ${\mathcal H}_0$-invariant  this form  is also
${\mathcal  H}_0$-invariant  and  so   it  follows  that  for  $[w]\in
H^d_c({\mathcal U})$, we have \[[v]\wedge  [ w] = I( [S_H(\Gamma _0)])
\wedge   [w].    \]   However   the   map   \[H^d_c({\mathcal   U})\ra
H^D_c(S_G)^{K_0}=j(H^D(\widehat  {X}))\]  given  by wedging  with  the
fundamental class $[S_H(\Gamma  _0)]$ is up to a  nonzero multiple the
integral  of  $[w]$ over  $\widehat  {Y}$ and  is  hence  equal to  $)
[\widehat   {Y}]\wedge  j^*([w]$  where   $j^*:H^*_c({\mathcal  U})\ra
H^*(\widehat  {X})$  is the  dual  of  the  restriction map  from  the
cohomology of $\widehat {X}$ to that of ${\mathcal U}$. By definition,
this means that \[ [S_H(\Gamma _0)]\wedge [w]= j([\widehat {Y}])\wedge
[w].\] This holds for all $[w] \in H^d_c({\mathcal U})$.  Hence
\[I([S_H(\Gamma _0)])=j([\widehat {Y}]).\] Here, under the isomorphism
of Franke mentioned in an earlier paragraph, the image of $j$ has been
identified with the  image of the restriction map  from the cohomology
of $\widehat {X}$ into that of ${\mathcal U}$.

This proves Theorem \ref{nonvan}.

\subsection{Remark} The foregoing proof is an adaptation of \cite{V}
where the  analogous result  for the compact  case is proved.   In the
compact  case,  the  fact  that  the  action  by  Hecke  operators  is
completely reducible  is used in a  crucial way. The  extension in the
present paper to  the non-compact case (i.e. when  $S(\Gamma )$ is not
compact) is achievable  only because of the result  of Franke that the
space  of Hecke algebra  invariants is  a {\bf  direct summand}  -as a
Hecke module- of $H^*(S_G)^{K_0}$ for the good maximal compact $K_0$.

\section{Applications}

In  this  section we  apply  the  criterion  of Theorem  \ref{nonvan}.
Before doing so  we note another result from  \cite{F} ( see \cite{F},
(7.1), equation (54)) on the kernel of the Borel map
\[j: H^*(\widehat {X}) \rightarrow H^*(S(\Gamma)).\]
We  have a  dual map  \[j^*:H^*_c(S(\Gamma))  \rightarrow H^*(\widehat
{X}).\] Franke's theorem says  that a class $[x]\in H^*(\widehat {X})$
lies   in   the   image   of  the   compactly   supported   cohomology
$H^*_c(S(\Gamma))$ if and  only if the restriction of  the class $[x]$
to  $\widehat  {X_L}$ vanishes  for  all  Levi  subgroups $L$  of  all
$\Q$-parabolic subgroups  of $G$.  Dually, this means  that the kernel
of the  Borel map  $j$ is  the orthogonal complement  of the  space of
vectors  $[v]\in  H^*(\widehat{X})$  whose  restriction  to  $\widehat
{X_L}$  vanishes for  al  the  Levis $L$  as  above.  Here  orthogonal
complement means  the following: the  cohomology algebra $H^*(\widehat
{X})$  comes equipped  with a  non-degenerate bilinear  form (Poincare
duality)  and  the  orthogonal  complement  is with  respect  to  this
bilinear form.

\subsection{Example} Let $G=R_{E/Q}(SL_n)$ with $E$ an imaginary
quadratic  extension  of $\Q$.   Then,  $\widehat  {X}=SU(n)$ and  its
cohomology is an exterior algebra
\[ \wedge=\wedge (e_3,e_5,\cdots, e_{2n-1})\]
on primitive generators $e_{2i-1}$ of degree $2i-1$.

Fix  $k\geq 1$ and  consider the  inclusion $SU(k)\subset  SU(n)$. The
cohomology algebra of $SU(k)$ is the exterior algebra
\[ \wedge (e_3, e_5,\cdots,e_{2k-1})\] and the restriction map in the
cohomology from  $SU(n)$ to  $SU(k)$ is given  by $e_i\mapsto  e_i$ if
$i\leq 2k-1$ and  $e_i\mapsto 0$ otherwise (this is  well known; e.g.,
see \cite{M-T}, Chapter (III), p.148, Theorem (6.5) (4)).

>From  the description  of the  Levi subgroups  (they are  of  the form
$R_{E/Q}(SL_{m_1}\times  \cdots  \times  SL_{m_k})$) it  follows  that
their  compact  duals  are   products  of  lower  dimensional  unitary
groups. Using this and the definition of the generators $e_{2i-1}$, it
can be proved  that a class $v\in \wedge$  restricts trivially to {\bf
all} these compact  duals if and only if it is  divisible by the class
$e_{2n-1}$.  It is  clear that the orthogonal complement  of the ideal
generated by  $e_{2n-1}$ is itself.  Hence by  Franke's theorem quoted
at  the beginning  of this  section, $j(v)=0$  if and  only if  $v$ is
divisible by $e_{2n-1}$.\\

Replace  $n$   by  $2n$  and  consider   the  embedding  $H=R_{E/Q}(Sp
_{2m})\subset G=R_{E/Q}(SL_{2m})$. The  compact dual $\widehat {Y}$ is
simply the group $Sp_{2n}$ whose cohomology is an exterior algebra
\[\wedge (e_3,e_7,\cdots, e_{4n-1}),\] on odd degree generators
$e_{2i-1}$  of  degree  $2i-1$  with  $i\geq 2$.   The  cohomology  of
$\widehat  {X}$ is  (note that  $n$ is  replaced by  $2n$) as  we said
before,
\[\wedge (e_3,e_5,\cdots,e_{4n-1}).\] Moreover, the restriction map
from $\widehat {X}$ to $\widehat {Y}$ takes $e_i$ to $e_i$ if $i=4j-1$
and  to $0$  otherwise.  This  can easily  be proved  by  looking more
closely at  the spectral  sequences used to  obtain the  cohomology of
$\widehat {X}$ and $\widehat  {Y}$ (see \cite{M-T}, p.119, Ch.  (III),
Theorem (3.10) (1) and (2)).

Consequently, the fundamental class is
\[ [\widehat {Y}]=e_5\wedge e_9\wedge \cdots \wedge e_{4n-3},\] and is
not divisible  by $e_{4n-1}$. This proves  that $j([\widehat {Y}])\neq
0$.  Therefore,  the modular symbol  $[S_H(\Gamma )]$ is  non-zero for
some congruence subgroup $\Gamma \subset G(\Q)$.

\medskip  Since  $[\widehat{Y}]$   is  not  divisible  by  $e_{4n-3}$,
Franke's result shows  that $j([\widehat{Y}])$ is not in  the image of
the cohomology with compact support. On the other hand the restriction
of  \[[\widehat {Y}]=e_5\wedge e_9\wedge  \cdots \wedge  e_{4n-3}\] to
the Levi subgroup $R_{E/Q}(SL_{n-1})$ is in the kernel of the image of
j for  $R_{E/Q}(SL_{n-1})$ since it  is divisible by  $e_{4n-3}$. Thus
$j([\widehat {Y}])$ is a ghost class.

\subsection{Example} Consider the embedding $SL_{2n+1}\subset
R_{E/\Q}(SL_{2n+1})$  of  $\Q$-groups,   where  $E$  is  an  imaginary
quadratic    extension    of    $\Q$.    The    embedding    $\widehat
{Y}=SU(2n+1)/SO(2n+1)\subset  \widehat  {X}=SU(2n+1)$  is  induced  by
$g\mapsto  gg^t$  from  $SU(2n+1)$  into itself.   The  cohomology  of
$SU(2n+1)/SO(2n+1)$  is (  see  \cite{F}, p.   35,  Proposition 7)  an
exterior algebra
\[\wedge (e_5, e_9,\cdots, e_{4n+1})\] and the restriction map from
\[H^*(\widehat {X})=\wedge (e_3,e_5,\cdots,e_{4n-1},e_{4n+1})\] to
$H^*(\widehat {Y})$  is given by  sending $e_{2i+1}$ to  $e_{2i+1}$ if
$i$ is odd and to $0$ otherwise.

Consequently, the fundamental class is
\[[\widehat {Y}]=e_3\wedge e_7\wedge \cdots \wedge e_{4n-1}\] and is
therefore not  divisible by $e_{4n+1}$. Hence  it does not  lie in the
kernel  of the  Borel map,  and by  Theorem \ref{nonvan},  the modular
symbol $[S_H(\Gamma )]$ is non-zero for a suitable $\Gamma $.\\

The rest of Theorem \ref{examples}  is proved in an entirely analogous
way. The  argument in  the previous example  shows that  $ j([\widehat
{Y}])=j(e_3\wedge e_7\wedge \cdots \wedge e_{4n-1})$ is a ghost class.

\medskip

\subsection{Remark} Suppose that $H$ is any linear
algebraic  group  over a  totally  real  number  field $F$  such  that
$H(F\otimes \R)=SL_{2n+1}(\R)^m$ with $m$  the degree of $F$ over $\Q$
(that is,  $H$ is  a $F$-form of  $SL_{2n+1}(\R)^m$).  Let $E/F$  be a
totally  imaginary quadratic extension  and set  $G=R_{E/F}(H)$. Then,
$[S_H(\Gamma  )]$   is  non-zero  for  some   congruence  subgroup  of
$G(E)$. This follows  by the same method as in  the case when $H=SL_n$
over $F$ (i.e. when $H$ is the standard $F$-form).\\

In the next two sections, we assume that both the symmetric spaces $X$
and $Y$ are  Hermitian symmetric and that the embedding  of $Y$ in $X$
is  holomorphic. Under  these assumptions,  one  can show  in a  large
number of cases that the associated modular symbol $[S_H(\Gamma )]$ is
non-zero for  some $\Gamma  $.  We will  discuss two  important cases,
when  $G$  is  the  symplectic  group $Sp_g$  and  the  unitary  group
$U(p,q)$.

\section{The Symplectic Group}

In  this  section,  $G=Sp_{2g}$  denotes  the  symplectic  group
defined and split over $\Q$. The associated symmetric space
is the Siegel upper  half space $X=Sp_{2g}(\R)/U(g)$. Denote by $\widehat
{X}$ its compact dual. Let $T$  be the group of diagonals in $K=U(g)$;
it is  a maximal torus  in $K$ and  in $G(\R)$. The cohomology  of the
classifying space $BT$  is generated by Chern classes  of line bundles
arising  from  characters of  $T$,  and  may  be identified  with  the
polynomial algebra $H^*(BT)=\C [\Lambda]$, where $X^*(T)$ is the group
of  characters on  $T$  and  $\Lambda \subset  X^*(T)$  is a  suitable
``positive   subset''  (in   the  case   of  $U(g)$),   consisting  of
non-negative  integral linear  combination  of the  characters of  $T$
occurring  in  the  standard   representation  of  $U(g)$  on  $\C^g$.
Clearly, the  ring $R=\C [\Lambda  ]$ is $\C [x_1,\cdots  ,x_g]$ where
$x_i$  is  (by a  mild  abuse  of notation)  the  Chern  class of  the
character $x_i$ (the  (i,i)-th entry) in the diagonal  torus $T$. Note
that   the   ring   $R$   is    a   module   for   the   Weyl   groups
$W_G=(\Z/2\Z)^g\times   S_g$   and    $W_K=S_g$   of   $G$   and   $K$
respectively. Here  $S_g$ is  the symmetric group  on the  $g$ letters
$x_1,x_2,  \cdots,  x_g$  and   $W_G$  is  a  semi-direct  product  of
$(\Z/2\Z)^g$ with  $S_g$; the  elements of $W_G$  act on the  $x_i$ by
$\pm   x_{\sigma  (i)}$,   with   $\sigma  \in   S_G$.   Let   $\sigma
_1,\cdots,\sigma _g$ denote the  elementary symmetric functions in the
variables $x_1,\cdots,x_g$.
\begin{lemma}\label{cohsp}The cohomology of $\widehat X$ is the ring
$\C[\sigma _1,\cdots,  \sigma _g]$ modulo  the ideal generated  by the
graded relation $\prod _i (1-x_i^2)=1$ ($i$ runs from $1$ to $g$).
\end{lemma}
\begin{proof} The cohomology of the compact dual $\widehat X$ may be
identified (\cite{M-T}) with the ring $R^{W_K}$ of $W_K$-invariants in
$R$ modulo the ideal generated by positive degree elements in the ring
$R^{W_G}$  of $W_G$-invariants.   The  above description  of the  Weyl
groups in question then implies the Lemma.
\end{proof}

Next, we determine the kernel of the Borel map, using the criterion of
Franke.   Let  $K=U(g)\subset  \widehat   {Sp_{2g}}$  be   the  natural
inclusion,  where  $G_u=\widehat {Sp_{2g}}$  is  the  compact form  of
$G=Sp_{2g}$.
\begin{lemma}\label{sptosu}
If there is {\bf any} embedding $i$ of the unitary group $U(g)$ in the
compact  form  $\widehat {Sp_{2g}}$  such  that  $K\cap i(SU(g))=O(g)$  the
orthogonal group and the restriction of the embedding $i$ to $O(g)$ is
identity,  then,  the  restriction   map  from  the  cohomology  (with
$\C$-coefficients) of $\widehat X$ into that of the subsymmetric space
$i(SU(g))/SO(g)$ is zero except in degree zero.
\end{lemma}
\begin{proof} To see this, first note that the cohomology of
$\widehat X$ is  generated by Chern classes of  the homogeneous vector
bundle arising from the {\bf standard} representation of $U(g)$ on $\C
^g$. Thus it is enough to  show that these Chern classes vanish on the
subsymmetric  space.   Now  the  restriction  $E$ of  this  bundle  to
$i(SU(g))/SO(g)$  is also  homogeneous, and  arises from  the standard
representation $\rho  $ of $SO(g)$.   Obviously, $\rho $ extends  to a
representation of  $i(SU(g))$, which shows that the  vector bundle $E$
admits a  trivialisation.  Hence all  its higher degree  Chern classes
are zero.  But  these classes are the restriction  of Chern classes of
the vector bundle  on $\widehat X$ with which  we started.  Therefore,
the restriction of the  cohomology of $\widehat X$ to $i(SU(g))/SO(g)$
is trivial.
\end{proof}

Let $P\subset G$ be a standard maximal parabolic subgroup defined over
$\Q$,  and  $L$  its  Levi  component.   One  may  identify  $L$  with
$GL_{g-k}\times  Sp_{2k}$ for  some $k\leq  g$. At  the level  of compact
duals,  one then  has  $\widehat {X_L}=SU(g-k)/SO(g-k)\times  \widehat
{Sp_{2k}}/U(k)$.   Thus  the  restriction  map  from  the  cohomology  of
$\widehat  {X}$  to  $\widehat  {X_L}$  factors  through  the  product
$\widehat  {Sp_{2g-2k}}/U(g-k)\times \widehat {Sp_{2k}}/U(k)$.   From Lemma
\ref{sptosu} (applied to  $Sp_{2g-2k}$) one sees that the  kernel of the
restriction map from $H^*(\widehat X)$ to $H^*(\widehat {X_L})$ is the
same as  the kernel of the  restriction map from  $H^*(\widehat X)$ to
$H^*(\widehat {Sp_{2k}}/U(k))$. \\

Let $J$ be the kernel of the restriction map $H^*(\widehat X)\ra \prod
H^*(\widehat  {X_L})$  where  the  product  runs  over  all  the  Levi
subgroups  $L$  of standard  maximal  parabolic  $\Q$- subgroups  $P$.
Therefore, we have proved that the kernel of the foregoing restriction
map  is  the  same  as  the  kernel of  the  map  $H^*(\widehat  X)\ra
H^*(\widehat {Sp_{2g-2}}/U(g-1))$.  Let $Z=\widehat {Sp_{2g-2}}/U(g-1)$,
and let  $\tau _1,\cdots,\tau _{g-1}$ denote  the elementary symmetric
functions   in  the   $g-1$-variables   $x_2,\cdots,x_g$.   By   Lemma
\ref{cohsp}, we have the identifications
\[H^*(\widehat X)=\C[\sigma _1,\cdots,\sigma
_g]/\prod _{i=1}^{i=g}(1-x_i^2)=1\] and
\[H^*(Z)=\C[\tau _1,\cdots,\tau
_{g-1}]/\prod _{i=2}^{i=g}(1-x_i^2)=1.\]

\begin{lemma}\label{sp} The kernel of the Borel map
$H^*(\widehat  X)\ra  H^*(S(\Gamma  ))$  is orthogonal  to  the  ideal
generated by $\sigma _g$.
\end{lemma}
\begin{proof}
By Franke,  the kernel  of the Borel  map is precisely  the orthogonal
complement (with respect to  Poincar\'{e} duality on $H^*\widehat X)$)
of  the  kernel of  the  map  $H^*(\widehat  X)\ra \prod  H^*(\widehat
{X_L})$.  By  the discussion  preceding the lemma,  the latter  is the
kernel  of  the   restriction  map  $H^*(\widehat  X)\ra  H^*(\widehat
{Sp_{2g-2}}/U(g-1))=H^*(Z)$.\\

Consider  the  map  $p$  from  the  ring  $\C[\sigma  _1,\cdots,\sigma
_g]/(\prod     _i(1-x_i^2)=1)$      into     the     ring     $\C[\tau
_1,\cdots,\tau_{g-1}]/(\prod   _j(1-x_j^2)-1)$  induced  by   the  map
$\C[x_,\cdots,   x_g]\ra  \C[x_2,\cdots,x_g]$,  with   $x_1\mapsto  0,
x_2\mapsto x_2,\cdots,x_g\mapsto x_g$.  Under the identifications made
just before Lemma \ref{sp}, $p$  is nothing but the restriction map of
the  previous  paragraph. The  kernel  of $p$  is  easily  seen to  be
generated by $\sigma _g$.
\end{proof}

In  the  proof of  Theorem  \ref{symplectic},  we  need the  following
lemmata.
\begin{lemma}\label{sympaux} If $k\leq g$, then in the ring $\C[\sigma
  _1,\cdots,\sigma  _g]/(\prod (1-x_i^2)-1)$,  the product  $\sigma _k
^2\sigma _{k+1}\cdots \sigma _g$ vanishes.
\end{lemma}
\begin{proof} By induction on $g$. Denote by $\beta $ the product
$\sigma  _k ^2\sigma _{k+1}\cdots  \sigma _g$.   Consider the  map $p$
into   $\C[\tau  _1,\cdots,\tau  _{g-1}]/(\prod   (1-x_j^2)-1)$.   The
element $\alpha  =\sigma _k^2\sigma _{k+1}\cdots  \sigma _{g-1}$ maps,
under  $p$, to the  element $\tau  _k^2\tau _{k+1}\cdots  \tau _{g-1}$
which is zero by  induction assumption.  Therefore, by Lemma \ref{sp},
the element $\alpha $ is divisible by $\sigma _g$. Now, $\beta =\alpha
\sigma _g$ and  is thus divisible by $\sigma  _g^2$.  However, $\sigma
_g^2$  is  the highest  degree  term  in  the graded  equation  $\prod
(1-x_i^2)=1$ and hence is zero. Therefore, $\beta $ also vanishes.
\end{proof}

\begin{lemma}\label{topsymp} The element
$\sigma _1\sigma _2\cdots\sigma _g$ is non-zero, i.e. it generates the
top degree (degree $g(g+1)$) cohomology of $\widehat X$.
\end{lemma}
\begin{proof} By induction on $g$. Since the only class in degree two
is  $\sigma _1$, it  is clear  that it  is the  {\bf Kahler  class} of
$\widehat X$. Hence, the top degree cohomology is generated by $\sigma
_1^{\frac {g(g+1)}{2}}$.  We will prove that $\sigma _1^{g(g+1)/2}$ is
a  multiple of  $\sigma  _1\cdots  \sigma _g$.   This  will prove  the
lemma. \\

Now, the  cohomology of $Z=\widehat {Sp_{2g-2}}/U(g)$  is the quotient
of that of $\widehat X$, by the ideal $\sigma _g$. Since the dimension
of $Z$ is $g(g-1)/2$,  it follows that $\sigma _1^{g(g+1)/2}$ vanishes
on  $Z$.   Thus,  in  the  ring $H^*(\widehat  X)$,  we  have  $\sigma
_1^{(g+1)g/2}=\sigma _g \psi$  where $\psi $ is a  degree $g(g-1)/2$ -
element.  Now, by induction assumption, the restriction of the element
$\psi$  to $Z$  is  a  multiple of  $\sigma  _1\cdots \sigma  _{g-1}$.
Consequently,   $\psi=\lambda   \sigma  _1\cdots\sigma   _{g-1}+\sigma
_g\phi$ for  some $\phi \in  H^*(\widehat X)$.  Multiplying  this last
equation by $\sigma  _g$ and noting that $\sigma  _g^2=0$, we see that
$\sigma _1^{g(g+1)/2}=\lambda \sigma  _1\cdots \sigma _g$, proving the
lemma.
\end{proof}
We will now begin the proof of Theorem \ref{symplectic}.

By the  criterion of Theorem \ref{nonvan},  it is enough  to show that
under the  Borel map, the image  of the compact  dual class $[\widehat
Y]$  is non-zero  , where  $\widehat  Y$ is  the compact  dual of  the
symmetric space of  $H$.  We will assume for simplicity  that $H$ is a
product of  {\bf two} symplectic  groups: put $a+b=g$ with  $a\geq b$,
and set $H=Sp_{2a}\times Sp_{2b}\subset Sp_{2g}$. The proof in the general case
is tedious, and we omit it.\\

Now,  the  cohomology   ring  $H^*(\widehat  {Sp_{2a}}/U(a))$  may  be
identified   with   $\C[x_1,\cdots,x_a]^{S_a}/(\prod  _i(1-x_i^2)=1)$,
where $S_a$ is the symmetric group on the $a$ letters $x_1,\cdots,x_a$
and $i$ runs from $1$ to  $a$. Let $\alpha _1,\cdots,\alpha _a$ be the
elementary symmetric functions in the variables $x_1,\cdots,x_a$.

Similarly, the  cohomology ring $H^*(\widehat  {Sp_{2b}}/U(b))$ may be
identified      with      $\C[x_{a+1},\cdots,x_{a+b}=x_g]^{S_b}/(\prod
_i(1-x_i^2)=1)$, where $S_b$ is the symmetric group on the $b$ letters
$x_{a+1},\cdots,x_{a+b}=x_g$ and $i$ runs  from $a+1$ to $a+b=g$.  Let
$\beta _1,\cdots,\beta  _b$ be  the elementary symmetric  functions in
the variables $x_{a+1},\cdots,x_g$.\\

The   top    degree   term   (of    degree   $2(\frac{a(a+1)}{2}+\frac
{b(b+1)}{2})$) in  the cohomology  group of the  product $H^*(\widehat
{Sp_{2a}}/U(a)\times  \widehat  {Sp_{2b}}/U(b))$  is easily  seen  (by
Lemma \ref{topsymp})  to be generated  by the element  $\gamma =\alpha
_1\cdots\alpha _a\otimes \beta _1\cdots\beta _b$.\\

The above  description of the cohomology  of the compact  duals of the
spaces  $Sp_{2g}/U(g)$, $Sp_{2a}/U(a)$  and  $Sp_{2b}/U(b)$ identifies
them as the rings  generated by certain elementary symmetric functions
modulo  the ideal  of  relations $\prod  _i(1-x_i^2)  =1$ for  certain
integers $i$ (where $i$ runs respectively  from $1$ to $g$, $1$ to $a$
and $a+1$ to $g=a+b$).  From this it  is easy to see that the image of
$\sigma _k$ in the cohomology of the product $\widehat {Sp_{2a}}\times
\widehat  {Sp_{2b}}/U(a)\times  U(b)$  is  the  sum  $\sum  _r  \alpha
_r\otimes \beta _{k-r}$ where $r$ runs from $0$ to $k$ (if $r\geq a+1$
then  $\alpha  _r=0$ by  convention).   In  particular,  the image  of
$\sigma _g$ is $\alpha _a\otimes  \beta _b$.  By induction on $k$, and
by  using  Lemma  \ref{sympaux}   applied  to  the  symplectic  groups
$Sp_{2a}$  and $Sp_{2b}$,  one  can  show that  the  image of  $\sigma
_g\sigma _{g-2}\cdots \sigma _{g-2k}$ is
\[\alpha _a\alpha _{a-1}\cdots \alpha _{a-r}
\otimes \beta  _b\beta _{b-1}\cdots \beta _{b-r}.\]  In particular, if
$k=b$ we  get that the  image of $\sigma _g\sigma  _{g-2}\cdots \sigma
_{g-2b}$ is the element
\[\alpha _a\alpha _{a-1}\cdots \alpha _{a-b}\otimes \beta _b\beta
_{b-1}\cdots \beta _1.\]

Thus, the top degree term $\gamma  $ (of the cohomology of the product
of  $\widehat {Sp_{2a}}/U(a)$  and $\widehat  {Sp_{2b}}/U(b)$)  of the
last but  one paragraph, lies  in the image  of the element  $\theta =
(\sigma _g\sigma _{g-2}\cdots \sigma _{g-2b})(\sigma _1\sigma _2\cdots
\sigma  _{a-b-1})$ under  the  restriction map  of  the cohomology  of
$\widehat   X$   to   the   cohomology  of   the   product   $\widehat
{Sp_{2a}}/U(a)\times \widehat {Sp_{2b}}/U(b)$.\\

By the definition of the  class $[\widehat Y]$ ($\widehat Y$ being the
compact  dual of  the symmetric  space associated  to  $Sp_{2a} \times
Sp_{2b}$),  this   implies  the  equality   $\theta  \wedge  [\widehat
Y]=\sigma _1\cdots  \sigma _g$, namely  a generator of the  top degree
(of degree  $2(\frac {g(g+1}{2})$)  cohomology of $\widehat  X$. Since
this wedge product is non-zero,  and $\theta $ is divisible by $\sigma
_g$, it follows that the  class $[\widehat Y]$ is {\bf not} orthogonal
to the ideal generated by $\sigma _g$.\\

Now, Lemma \ref{sp}  implies that $[\widehat Y]$ is  not in the kernel
of the Borel map. This proves Theorem \ref{symplectic}.

\section{The Unitary Group}

In this section,  $G=U(p,q)$ will denote the unitary  group in $n=p+q$
variables, with  $1\leq p\leq  q$.  The $\Q$  structure is  defined as
follows.  Let $V$ be an $n$-dimensional vector space over an imaginary
quadratic extension  $E$ over $\Q$. Consider  the $E$-valued Hermitian
form in $n$ variables given by
\[h(v,v)=\sum _{i=1}^{i=p}\mid z_i\mid ^2-\sum _{i=p+1}^{i=n}\mid
z_i\mid ^2\] with $v=(z_1,\cdots,z_n)\in V$. The group preserving this
Hermitian  form is a  $Q$-algebraic group  and the  group of  its real
points is the group $U(p,q)$.\\

The  group  $K=U(p)\times  U(q)$  is  a maximal  compact  subgroup  of
$U(p,q)$ and the  group $T$ of diagonals in $K$ is  a maximal torus of
$G$ and $K$.  As in section (4.1), if $\widehat X$ is the compact dual
of  $G/K$  ($\widehat   X$  is  the  Grassmanian  of   $p$  planes  in
$n$-dimensional  complex   vector  space),  then   the  cohomology  of
$\widehat X$ may be identified with the quotient ring (see \cite{M-T})
\[\frac {\C[x_1,\cdots,x_p;~y_i,\cdots,y_q]^{S_p\times S_q}}
{\prod (1+x_i)\prod(1+y_j)=1}\] In  this equality, $S_p$ (resp, $S_q$)
is the  permutation group of  the $p$ letters  $x_1,\cdots,x_p$ (resp.
the $q$  letters $y_1,\cdots,y_q$).  The superscript  denotes the ring
of invariants  under the product group $S_p\times  S_q$.  The variable
$i$  (resp, $j$)  runs  from $1$  to  $p$ (resp.   $1$  to $q$).   The
equation $\prod (1+x_i)\prod (1+y_j)=1$ is a graded equation.\\

Let $\sigma  _1,\cdots,,\sigma _p$  (resp. $\tau _1,\cdots,  \tau _q$)
denote   the   elementary  symmetric   functions   in  the   variables
$x_1,\cdots,x_p$  (resp. $y_1,\cdots,y_q$).   Then, the  cohomology of
$\widehat X$ is the ring
\[\frac {\C[\sigma _1,\cdots, \sigma _p;~\tau _1,\cdots,\tau _q]}
{(1+\sigma  _1+\cdots  +\sigma _p)(1+\tau  _1  +\cdots +\tau  _q)=1}\]
where as before, the equation involving $\sigma $'s and $\tau $'s is a
graded one.
\begin{lemma}\label{unitborel} With the foregoing notation, the kernel
of the  Borel map for $G=U(p,q)$  is the orthogonal  complement of the
ideal generated  by $\sigma _p$ and  $\tau _q$ in  the cohomology ring
$H^*(\widehat X)$.
\end{lemma}
\begin{proof} As $p\leq q$, the $\Q$-rank of the group $G=U(p,q)$ is
  $p$.  The  Hermitian form $h$  defining $G$ is  a direct sum  of $p$
``hyperbolic'' (i.e.  isotropic  over $\Q$) Hermitian forms $h_1\oplus
\cdots \oplus  h_p$ and  an anisotropic hermitian  form $h_0$.   It is
then  easy to  see that  a Levi  subgroup is  of the  form $GL_k\times
U(p-k,q-k)$ for some $k\leq p$.  Thus, $\widehat {X_L}$ is the product
of  the symmetric  spaces  $SU(k)/SO(k)$ and  $U(p+q-2k)/(U(p-k)\times
U(q-k))$.\\

An  argument involving  Chern  classes of  homogeneous vector  bundles
similar to that in the proof  of Lemma \ref{sptosu} shows that all the
(positive degree) cohomology classes on $U(2k)/U(k)\times U(k)$ vanish
on $i(SU(k))/SO(k)$, for {\bf any} embedding $i$ of $SU(k)$ in $U(2k)$
such that the  intersection $i(SU(k))\cap (U(k)\times U(k))=SO(k)$ and
$i$  is the identity  on $SO(k)$.   In particular,  it shows  (cf. the
paragraph preceding Lemma \ref{sp}) that the kernel of the map
\[H^*(\widehat X)\ra \prod H^*(\widehat X_L)\]
(where the product is over all the standard Levi subgroups as before),
is the kernel of the map
\[H^*(\widehat X)\ra H^*(U(p+q-2)/(U(p-1)\times U(q-1))).\]
The latter  kernel may easily  be shown to  be the ideal  generated by
$\sigma _p$ and $\tau _q$. \\

By Franke,  the kernel  of the Borel  map is precisely  the orthogonal
complement of this ideal. Hence the Lemma.
\end{proof}

\begin{lemma}\label{unitl} With the previous notation, the element
$\tau _q^p$ generates  the top degree cohomology (in  degree $2pq$) of
$\widehat X$.
\end{lemma}
\begin{proof} By induction on $p$.

Consider the  restriction map from $H^*(\widehat  X)$ to $H^*(\widehat
Y_1)$, where $\widehat Y_1$ is the compact dual of the symmetric space
associated  to $U(p-1,q)\subset  U(p,q)$. The  kernel of  this  map is
generated by  $\sigma _p$,  as may be  easily seen. In  particular, if
$R_{pq-p-q}$  be  the graded  piece  of  the  cohomology ring  $R$  of
$\widehat X$ (since $H^{(p-1)q}(\widehat  Y_1)$ is one dimensional) it
follows  that  $\sigma  _p   R_{pq-p-q}$  is  of  codimension  one  in
$H^{(p-1)q}  (\widehat X)$.   But,  wedging by  $\tau  _q$ kills  this
subspace  (since, by the  graded relation,  $\tau _q\sigma  _p=0$). By
definition,  wedging by the  cycle class  $[\widehat Y_1]$  also kills
this subspace.   Therefore, $[\widehat Y_1]=\lambda \tau  _q$ for some
non-zero scalar $\lambda$.\\

By  induction, the restriction  of $\tau  _q^{p-1}$ generates  the top
degree cohomology of  $\widehat Y_1$.  By the definition  of the cycle
class  $[\widehat Y_1]$,  $\tau _q^{p-1}[\widehat  Y_1]$  is non-zero,
i.e.  generates  the top  degree cohomology of  $\widehat X$.   By the
previous  paragraph, this element  is, up  to non-zero  scalars, $\tau
_q^{p}$.
\end{proof}

\begin{remark}\label{unitr}Note that since $\sigma _1$ is the Kahler
class on  $\widehat X$  and restricts to  a Kahler class  on $\widehat
Y_1$,  we  have also  proved  that  $\tau _q(\sigma  _1)^{(p-1)q}=\tau
_q^p\neq 0$.
\end{remark}

We now begin the proof of Theorem \ref{unitary}.

\subsection{Part 1 of Theorem \ref{unitary}.} The Hermitian space
$V$ may be split into a  direct sum of Hermitian spaces $V_i$ and $W$,
where, on each $V_i$, the form  $h$ is of type $(p_i,q_i)$; on $W$, it
is  of type  $(0,q-\sum q_i)$.   Then $\sum  p_i=p$ and  $\sum q_i\leq
q$. We thus get an embedding of $H=\prod U(p_i,q_i)$ in $U(p,q)$.\\

Let $\xi _1,\cdots,  \xi _l$ (resp.  $\omega _1,\cdots  , \omega _l$ )
denote  the analogues of  $\tau _q$  (resp.  of  $\sigma _1$)  for the
groups $U(p_1,q_1),\cdots, U(p_l,q_l)$  and $Z=Z_1\times\cdots Z_l$ be
the products of the compact  symmetric spaces associted to the product
group. The  restriction of  $\tau _q$ to  this product variety  $Z$ is
clearly the tensor product $\xi _1\otimes \cdots \otimes \xi _l$. Here
we use the  fact that $\sum p_i=p$. An easy  computation shows that if
$d$  (resp $d_i$)  is the  dimension of  $Z$ (resp.   $Z_i$)  then the
restriction  of the  element  $\tau  _q \sigma  _1^d$  to the  product
variety    $Z$     is    the    tensor     product    $\xi    _1\omega
_1^{d_1}\otimes\cdots\otimes   \xi  _l\omega  _l^{d_l}$.    By  Remark
\ref{unitr} applied  to $U(p_i,q_i)$  for each $i$,  we get  that this
tensor  product element  is the  top degree  cohomology class  of $Z$.
Thus, $\tau  _q\sigma _1^{d}\wedge  [Z]\neq 0$.  Therefore,  the cycle
class  of the  compact dual  associated to  $\prod U(p_i,q_i)$  is not
orthogonal to  an element  of $\tau _qR$.   In particular  (from Lemma
\ref{unitborel}),  $[Z]$ does  not  lie  in the  kernel  of the  Borel
map. This  proves , by  Theorem \ref{nonvan}, that the  modular symbol
$[S_H(\Gamma )]$ is non-zero.

\subsection{Part (2) of Theorem \ref{unitary}.} The cohomology of
the  compact dual  $\widehat Y$  associated  to $Sp_{2g}$  is (see  Lemma
\ref{cohsp})
\[\frac {\C[\sigma _1,\cdots, \sigma _g]}{\prod (1-x_i^2)=1}.\]
The cohomology of $\widehat X$ for the group $U(g,g)$ is
\[\frac {\C[\sigma _1,\cdots,\sigma _g; \tau _1,\cdots,\tau _g]}
{\prod (1+x_i)(1+y_i)=1}\] The restriction  map from the cohomology of
$\widehat X$ to  that of $\widehat Y$, is  induced by $x_i\mapsto x_i$
and  $y_i\mapsto  -x_i$.   Therefore,  the top  degree  class  $\sigma
_1\cdots\sigma  _g$ of $\widehat  Y$ is  in the  image of  the product
$\tau _1\cdots\tau _g$. Hence,  $[\widehat Y]\wedge \tau _1\cdots \tau
_g$ generates  the top degree class  of $\widehat X$.  This means that
the  cycle  class  $[\widehat  Y]$  is not  orthogonal  to  the  ideal
generated by  $\tau _g$.   Thus, by Lemma  \ref{unitborel}, $[\widehat
Y]$  is  not in  the  kernel  of the  Borel  map,  whence, by  Theorem
\ref{nonvan}, the modular symbol $[S_H(\Gamma )]$ is non-zero.  \\

\subsection{Proof of Theorem \ref{infinite}} The Hermitian space
$(E,h)$  is  such   that  for  the  associated  group   $G$,  we  have
$G(\R)=SU(1,q)$. We first prove the following lemma.

\begin{lemma} \label{unitary(1)}
If  $\alpha $  is a  non-zero  holomorphic $1$-form  on $S(\Gamma  )$,
$\omega $ is  the Kahler class on $\widehat X$ and  $j$ the Borel map,
then the cup-product $\alpha \wedge j(\omega )$ is non-zero.
\end{lemma}
\begin{proof} It is enough to prove (since $\omega $ is $G(\A_f)$-
invariant), that  for some  $g\in G(\A_f)$, the  cup-product $g(\alpha
)\wedge  \omega $ is  non-zero. Consider  the {\bf  three} dimensional
hermitian subspace  $(F,h\mid _F)$ of the Hermitian  space $(V,h)$, (a
vector space  over the  field $E$) where  the Hermitian form  does not
represent a zero.  By weak  approximation, this is possible (since one
may locate a  three dimensional subspace over $p$-adic  field $E_v$ of
the  hermitian  vector space  $V\otimes  E_v,h\otimes  E_v$ where  the
Hermitian form does not represent a zero). \\

This  is   an  anisotropic  subspace,  whence   the  associated  group
$H_0=SU(F,h)$ is anisotropic over  $\Q$ and is isomorphic to $SU(1,2)$
over $\R$.  Thus, the  subvariety $S_{H_0}(\Gamma )$ is {\bf compact}.
If $\alpha $ is a  non-zero holomorphic $1$-form on $S(\Gamma )$ class
on  $\widehat X$, then  there exists  a $g\in  G(\A_f)$ such  that the
restriction of the {\bf form}  $g^*(\alpha )$ to $S_{H_0}(\Gamma )$ is
non-zero.  By  replacing $\alpha  $ by $g^*(\alpha  )$, we  may assume
that $g=1$.  However, since  $S_{H_0}(\Gamma )$ is compact, this means
that  the  restriction of  $\alpha  \wedge  \overline {\alpha  }\wedge
j(\omega )$ to $S_{H_0}(\Gamma )$ is non-zero.  In particular, $\alpha
\wedge j(\omega )\neq 0$.\\

In  the  above,  $\overline  {\alpha}$  denotes  the  anti-holomorphic
$1$-form  which  is the  complex  conjugate  of  the holomorphic  form
$\alpha$.   The  complex  conjugation   is  on  the  cohomology  group
$H^1(S(\Gamma ),\C)=H^1(S(\Gamma ),\R)\otimes \C$.
\end{proof}

\begin{proof} {\it of Theorem  \ref{infinite}.} We will argue by contradiction.
 Suppose that for
$H=SU(1,q-1)$ the cycle class $[S_H(\Gamma )]$ is always (i.e. for
{\bf every} $\Gamma $) $G(\A_f)$-invariant. By Theorem
\ref{nonvan},  this class  is then equal  to $j([\widehat Y])$.
However, $H^2(\widehat X)$  is one dimensional ($\widehat X={\bf
P}^q$, the complex projective $q$-space). Hence $[\widehat
Y]=\omega $ the Kahler class.

If $\alpha $ is a non-zero holomorphic form on $S(\Gamma )$ (such
forms exist  by \cite{K}), then $\alpha \wedge  j([\widehat
Y])\neq 0$ by  Lemma \ref{unitary(1)}.   Thus, the  restriction of
$\alpha  $ to $S_H(\Gamma )$ does not vanish, for {\bf every}
$\Gamma $. Thus,  at the level  of $H^1$  the restriction  map
from  $SU(1,q)$ to $SU(1,q-1)$ is injective for {\bf every}
$\Gamma $. In particular, the holomorphic cohomology  classes of
degree one  associated to $SU(n,1)$ restrict injectively  to those
on $SU(n-1,1)$.

We now recall the criterion of \cite{V2}, for the span of
Hecke-translates of $[S_H(\Gamma]$ to be infinite dimensional. In the
statement of the following theorem, the map $Res$ refers to an ``Oda
style'' restriction map.

\begin{theorem}(see \cite{V2}, Theorem 1). Suppose that $G$ is almost
  $\Q$-simple and that

1. The centralizer $Z_G(H)\cap K$ is not contained in the center
   $Z(G)$ of $G$.

2. For some integer $m\leq d=dim (S_H(\Gamma ))$ (dimension as a
   complex manifold), the restriction map
\[Res: H^{m,0}(S(\Gamma ))
   \ra \prod _{g\in G(\Q)} H^{m,0}(S_H(g\Gamma g^{-1}))\]
is non-zero.

Then, there exists a congruence subgroup $\Gamma '$ of $\Gamma$ such
that the cycle class $[S_H(\Gamma ')]$ is not $G(\A_f)$ invariant.
\end{theorem}

For the group $G$ with $G(\R)=SU(n,1)$ up to compact factors there
exist elements in $G(\Q)$  which centralize $SU(n-1,1)$ but do not
lie in  the center  of $G$, since the centralizer of $SU(n-1,1)$
is the group $U(1)$ (all  these viewed, by restriction of scalars,
as groups over $\Q$).

Thus, the conditions of Theorem 1 of \cite{V2} are satisfied and
so by Theorem 1 of \cite{V2}, there exists a $\Gamma$ so that  the
cycle class
$[S_H(\Gamma )]$ is not $G(\A_f)$-invariant, which contradicts our assumption.
\end{proof}

\noindent{\bf Acknowledgement.} This work was begun when both the
authors were visiting MPI, Bonn  in May-June of 2001; the hospitality
of MPI is gratefully acknowledged.  The authors also thank Arvind Nair
for very helpful conversations. Birgit Speh was partially supported by
the NSF grant DMS-007056. \\

\end{document}